\documentclass{amsart} 
\usepackage{epsf} 
\usepackage{amssymb,latexsym, amsmath, amscd, array, graphicx} 
\def\hepsffile{\leavevmode\epsffile}

\swapnumbers 
\numberwithin{equation}{section} 

\theoremstyle{plain} 
 
\newtheorem{thm}{Theorem}[section]

\theoremstyle{definition} 
\newtheorem{defin}[thm]{Definition}

\newtheorem{rem}[thm]{Remark} 
 
\newtheorem{ex}[thm]{Example}

\def\id{\protect\operatorname{id}} 
\def\Im{\protect\operatorname{Im}}

\def\span{\protect\operatorname{span}}

\def\pr{\protect\operatorname{pr}}

\def\exp{\protect\operatorname{exp}}

\def\pt{\protect\operatorname{pt}}

\def\R{{\mathbb R}}

\def\ss{{\sf X}}

\def\1{\hbox{\rm\rlap {1}\hskip.03in{\rom I}}} 
\def\Bbbone{{\rm1\mathchoice{\kern-0.25em}{\kern-0.25em} 
{\kern-0.2em}{\kern-0.2em}I}} 

\def\pp{\medskip{\parindent 0pt \it Proof.\ }} 
\def\wt{\widetilde} 
\def\wh{\widehat} 
 
\def\ov{\overline} 
 
\long\def\forget#1\forgotten{} %

\begin{document} 
\date{February 10, 2007} 
\leftline{ } 
\centerline{ } 
\title[Mapped Null Hypersurfaces and Legendrian Maps] 
{Mapped Null Hypersurfaces and Legendrian Maps} 
\author[V.~Chernov (Tchernov)]{Vladimir V. 
Chernov 
(Tchernov)} 
\address{V. Chernov, Department of Mathematics, 
6188 Kemeny Hall, Dartmouth College, Hanover NH 03755, 
USA} 
\email{Vladimir.Chernov@dartmouth.edu} 

\begin{abstract}
For an $(m+1)$-dimensional space-time $(\ss^{m+1}, g),$  define a mapped null hypersurface to be a smooth map $\nu:\mathcal N^{m}\to \ss^{m+1}$ (that is not necessarily an immersion) such that there exists a smooth field of null lines along $\nu$ that are both tangent and $g$-orthogonal to $\nu.$  

We study relations between mapped null hypersurfaces and Legendrian maps to the spherical cotangent bundle $ST^*\mathcal M$ of an immersed spacelike hypersurface $\mu:\mathcal M^m\to \ss^{m+1}.$  We show that a Legendrian map $\wt \lambda: \mathcal L^{m-1}\to (ST^*\mathcal M)^{2m-1}$ defines a mapped null hypersurface in $\ss.$ On the other hand, the intersection of a mapped null hypersurface $\nu:\mathcal N^m\to \ss^{m+1}$ with an immersed spacelike hypersurface $\mu':\mathcal M'^m\to \ss^{m+1}$ defines a Legendrian map to the spherical cotangent bundle $ST^*\mathcal M'.$ This map is a Legendrian immersion if $\nu$ came from a Legendrian immersion to $ST^*\mathcal M$ for some immersed spacelike hypersurface $\mu:\mathcal M^m\to \ss^{m+1}.$
\end{abstract}

\thanks{2000 Mathematics Subject Classification: Primary 53C50; Secondary 57R17} 

\keywords{Lorentz manifold, hypersurface, null surface, null cone, contact manifold, Legendrian manifold} 
\maketitle

We work in the $C^{\infty}$ category, and the word ``smooth'' means  
$C^{\infty}$. The manifolds in this work are assumed to be smooth without boundary. They are not assumed to be oriented, or connected, or compact unless the opposite is explicitly stated. In this work $(\ss^{m+1}, g)$ is an $(m+1)$-dimensional Lorentzian manifold that is not assumed to be geodesically complete. 

A ``vector field'' on a manifold $Y$ is a smooth section of the tangent bundle $\tau_{Y}:TY\to Y$, and a ``vector field along a map'' $\phi:Y_1\to Y_2$ of one manifold to another is a smooth map $\Phi:Y_1\to TY_2$ such that $\phi=\tau_{Y_2}\circ \Phi.$ Covector fields and line fields on a manifold and along a map $\phi$ are defined in a similar way.

\section{Preliminaries}
{\it Let us recall some basic Lorentz geometry facts.\/}
Put $\Xi$ to be the space of vector fields on $\ss.$  There exists a unique connection $\nabla^g$ on $\ss$ that satisfies the following metric compatibility and torsion free conditions:

\begin{equation}
\begin{split}
\xi_1 g(\xi_2, \xi_3)=g(\nabla^g_{\xi_1}\xi_2, \xi_3)+g(\xi_2, \nabla^g_{\xi_1}\xi_3), \\ [\xi_1, \xi_2]=\nabla^g_{\xi_1}\xi_2-\nabla^g_{\xi_2}\xi_1,
\end{split}
\end{equation}
for all $\xi_1, \xi_2, \xi_3\in \Xi;$ see~\cite[page 22]{BeemEhrlichEasley}. This connection is called a {\it Levi-Civita connection.\/} When no confusion arises we will write $\nabla$ instead of $\nabla^g.$

A {\em geodesic\/} $c:(a,b)\to (\ss, g)$ is a  
smooth curve satisfying $\nabla_{  c'}  c'=0,$ for all of its points.
Similarly to the Riemannian case, one uses geodesics to define the exponential 
$\exp_x:T_x\ss\to \ss.$ Note that $\exp_x$ is defined on the star-convex, with respect to ${\bf 0},$ domain of $T_x\ss,$ rather than on the whole $T_x\ss.$ 
There is an open neighborhood $V_x\subset T_x\ss$ of ${\bf 0}$ such that $\exp_x|_{V_x}$ is a diffeomorphism. The open set $U_x=\exp_x(V_x)$ is called a {\it normal neighborhood\/} of $x.$ 
A neighborhood is {\it geodesically convex\/} if any two of its points can be joined by a unique geodesic arc inside of it. The result of Whitehead~\cite{Whitehead},~\cite{Whiteheadaddendum},~\cite[Section 5, Proposition 7]{ONeill} is that every point in a semi-Riemannian, and hence Lorentzian manifold has a geodesically convex normal   neighborhood. 
A {\it simple region\/} is a geodesically convex normal neighborhood with compact closure whose boundary is diffeomorphic to $S^{m}.$

A nonzero vector $v\in T\ss$ is called {\em spacelike,  
non-spacelike, null $($lightlike$)$, {\rm or} timelike\/} if $g(v,v)$ is 
positive, non-positive,  zero, or negative, respectively. A piecewise smooth curve is called {\em spacelike,  non-spacelike, null, {\rm or} timelike\/} if all of its velocity vectors are respectively spacelike, non-spacelike, null, or timelike. For a point $x$ in a Lorentz $(\ss, g)$ the set of all nonspacelike vectors in $T_x\ss$ consists of two connected components that are hemicones. A continuous with respect to $x\in \ss$ choice of one of the two hemicones is called the {\it time orientation\/} of $\ss.$ The nonspacelike vectors from these chosen hemicones are called {\it future pointing.\/} A time oriented $(\ss^{m+1}, g)$ is called a {\it space-time.\/}

An immersion  
$\kappa:\mathcal K^k\to \ss^{m+1}$ of a $k$-manifold is said to be an immersed spacelike or timelike submanifold if the pull back of $g$ to $T\mathcal K$ is respectively a Riemannian or a Lorentzian metric.
An immersion (respectively an embedding) $i:\mathcal H^m\to \ss^{m+1}$ of an $m$-manifold is called an {\it immersed (respectively  an embedded) hypersurface.\/} 
An immersed hypersurface is called an {\it immersed null hypersurface\/} if for every $h\in \mathcal H$ the pull back of $g$ is degenerate on $T_h\mathcal H.$ Similarly one defines {\it embedded null hypersurfaces\/} and immersed and embedded 
{\it spacelike and timelike hypersurfaces.\/}

An immersed (or an embedded) hypersurface $i:\mathcal H^m\to \ss^{m+1}$ can be canonically equipped with a line field $L_{h}\subset T_{i(h)}\ss, h\in \mathcal H,$ along $i$ such that for every $h\in \mathcal H$ the line $L_h$ is $g$-orthogonal to $i_*(T_h\mathcal H)\subset T_{i(h)}\ss.$ It is easy to verify that an immersed hypersurface is spacelike, timelike or null if and only if for every $h\in \mathcal H$ the nonzero vectors in $L_h$ are respectively timelike, spacelike or null. Since the Lorentz metric is non-degenerate, for an immersed null surface the line field $L_h$ is tangent to $i(\mathcal H)$, i.e. $L_h\subset i_*(T_h\mathcal H)$ for all $h\in \mathcal H.$ This observation motivates the following definition.

\begin{defin}[mapped null hypersurface]
A smooth map $\nu:\mathcal N^{m}\to \ss^{m+1}$ of an $m$-manifold is called {\it a mapped null hypersurface\/} if there exists a smooth (non-oriented) line field $L_n\subset T_{\nu(n)}\ss, n\in \mathcal N,$ along $\nu$ such that for every $n\in \mathcal N$ the nonzero vectors in $L_n$ are null, $L_n\subset \nu_*(T_n\mathcal N),$ and $L_n$ is $g$-orthogonal to $\nu_*(T_n\mathcal N)\subset T_{\nu(n)}\ss.$

Two null vectors are orthogonal if and only if one of them is a multiple of the other. Hence the line field $L_n$ in the above definition is completely determined by the map $\nu.$

Every immersed null hypersurface is a mapped null hypersurface. However mapped null hypersurfaces can be quite singular. For example if $Y^{m-1}$ is an $(m-1)$-manifold and $\gamma:\R\to \ss$ is a curve such that $\gamma'(t)$ is null and nonzero for all $t,$  then the composition of $\gamma$ and of the projection $Y^{m-1}\times \R\to \R$ gives a mapped null hypersurface $Y^{m-1}\times \R\to \ss^{m+1}.$
\end{defin}

{\it Let us recall some basic contact geometry facts.\/} 
Let $Q^{2k-1}$ be  
a  smooth $(2k-1)$-dimensional manifold  
equipped with a smooth (non-oriented) hyperplane field $\zeta=\{\zeta^{2k-2}_q\subset  
T_qQ^{2k-1}\bigm|q\in Q\}.$ This hyperplane field is called {\em a contact  
structure,\/}  
if it can be locally presented as the kernel of a $1$-form $\alpha$ with nowhere zero $\alpha\wedge (d\alpha)^{k-1}.$  
 
An immersion (respectively an embedding) $i: \mathcal L^{k-1}\to Q^{2k-1}$ of a
$(k-1)$-dimensional  
manifold $\mathcal L^{k-1}$ into a contact manifold 
$(Q^{2k-1},\zeta)$  
is called a {\em Legendrian immersion (respectively a Legendrian embedding),\/} 
if $i_*(T_l\mathcal L)\subset  
\zeta_{i(l)},$ for all $l\in \mathcal L.$  

\begin{defin}[Legendrian map]
We say that a smooth map $\wt \lambda:\mathcal L^{k-1}\to Q^{2k-1}$ to a contact $(Q^{2k-1}, \zeta)$ is a {\it Legendrian map\/} if $\wt \lambda_*(T_l\mathcal L)\subset  
\zeta_{\wt \lambda(l)},$ for all $l\in \mathcal L.$ Every Legendrian immersion  is a Legendrian map. However a Legendrian map can be quite singular and the trivial map $\mathcal L^{k-1}\to \pt\in Q^{2k-1}$ is a Legendrian map.
\end{defin}
 
\begin{ex}[The natural contact structure on $ST^*\mathcal M$]\label{contactST*M}  
For a smooth manifold $\mathcal M^k,$ put $ST^*\mathcal M$ to be the spherical cotangent bundle, i.e. the quotient of $T^*\mathcal M$ minus the zero section by the action of the group $\R^+$ of positive real numbers under multiplication. Put $\pr=\pr_M:ST^*\mathcal M\to \mathcal M$ to be the corresponding $S^{k-1}$-bundle map. A point $p\in ST^*\mathcal M$ is the equivalence class of nonzero linear functionals on $T_{\pr p}\mathcal M.$ Two functionals are equivalent if and only if their kernels are equal and the half spaces of $T_{\pr p}\mathcal M$ where the functionals are positive are equal.
Thus $p$ is completely determined by the hyperplane $\ker 
p\subset  T_{\pr (p)}\mathcal M$ together with the halfspace of $T_{\pr p}\mathcal M\setminus \ker p$ where the functionals are positive. 

The natural contact structure  
$$ 
\zeta=\{\zeta_p^{2k-2}\subset T_p(ST^*\mathcal M)^{2k-1}, p\in  
ST^*\mathcal M\} 
$$  
is given by $\zeta_p=(\pr_*)^{-1}(\ker p)$. 

A map $\wt \lambda:\mathcal L^{k-1}\to (ST^*\mathcal M)^{2k-1}$ can be described as the pair consisting of the smooth map $ \lambda=\pr \circ \wt \lambda$ and a smooth nowhere zero covector field $\theta_l\in T^*_{ \lambda(l)}\mathcal M, l\in \mathcal L,$ along $ \lambda$ such that for every $l\in \mathcal L$ the equivalence class of $\theta_l$ is $\wt \lambda(l).$
The covector field $\theta_l$ is defined uniquely up to a multiplication by a positive smooth function $\mathcal L\to \R.$
 
Clearly $\wt \lambda$ is a Legendrian map if and only if $ \lambda_*(T_l\mathcal L)\subset \ker \theta_l,$ for all $l\in \mathcal L.$

If $\mathcal M$ is equipped with a  
Riemannian or Lorentzian metric $h$, then we can identify the tangent and the cotangent  
bundles of $\mathcal M$ and we can identify the spherical tangent and the spherical cotangent bundles $\pr:ST\mathcal M\to \mathcal M$ and $\pr:ST^*\mathcal M\to \mathcal M.$ Thus a smooth map $\wt \lambda:\mathcal L\to ST\mathcal M=ST^*\mathcal M$ can be described as the pair consisting of the smooth map $ \lambda=\pr \circ \wt \lambda$ and a smooth nowhere zero vector field $X_l\subset T_{ \lambda(l)}\mathcal M, l\in \mathcal L,$
along $ \lambda$ such that for every $l\in \mathcal L$ the equivalence class of $X_l$ is $\wt \lambda(l).$  Clearly $\wt \lambda$ is a Legendrian map if and only if $X(l)$ is $h$-orthogonal to $ \lambda_*(T_l\mathcal L),$ for all $l\in \mathcal L.$

Now let $(\ss^{m+1}, g)$ be a space-time and let $\nu:\mathcal N^m\to \ss^{m+1}$ be a mapped null hypersurface. Let $L_n\in T_{\nu(n)}\ss, n\in \mathcal N,$ be the unique smooth nonoriented line field along $\nu$ from the definition of the mapped null hypersurface. Since $(\ss, g)$ is time oriented, we can orient the null lines $L_n$ in the direction of the future. This oriented line field defines the Legendrian map $\wt \nu:\mathcal N\to ST\ss$ such that $\pr_{\ss}\circ \wt \nu=\nu.$ 
\end{ex}

\section{From Legendrian maps to mapped null hypersurfaces}
Let $(\ss^{m+1}, g)$ be a space-time,  
let $\mu:\mathcal M^{m}\to \ss^{m+1}$ be an immersed spacelike hypersurface, and let $\ov g$ be the induced Riemannian metric on $\mathcal M.$
Let $\wt \lambda:\mathcal L^{m-1}\to ST^*M=STM$ be a Legendrian map that is described
by the pair $ \lambda=\pr \circ \lambda$ and the smooth unit length vector field $X_{l}\in T_{ \lambda(l)}\mathcal M, l\in \mathcal L,$ along $ \lambda.$ 

Since the immersed hypersurface $\mu$ is spacelike, for each $l\in \mathcal L$ the space $T_{\mu\circ  \lambda (l)}\ss$ splits as the direct sum of $\mu_*(T_{ \lambda(l)}\mathcal M)$ and its one-dimensional $g$-orthogonal compliment $(\mu_*(T_{ \lambda(l)}\mathcal M))^{\perp}$ consisting of timelike vectors. Thus for each $l\in \mathcal L,$ there exists the unique future pointing null vector $N_l=(N_l^s, N_l^t)\in \mu_*(T_{ \lambda(l)}\mathcal M)\oplus (\mu_*(T_{ \lambda(l)}\mathcal M))^{\perp}=T_{\mu\circ  \lambda(l)}\ss$ such that $N_l^s=\mu_*(X_l).$ 
Put $\gamma_l(t)$ to be the maximal null geodesic such that $\gamma'_l(0)=N_l.$

We get the map from a subset of $\mathcal L\times \R$ to $\ss$ defined as 
$(l, t)\to \gamma_l(t),$ for $l\in \mathcal L$ and $t$ in the domain of the null geodesic $\gamma_l.$ Since each point of $\ss$ has a  geodesically convex normal neighborhood, the above map is defined on an open neighborhood of $\mathcal L\times 0\subset \mathcal L\times \R.$ Put $\mathcal N\subset \mathcal L^{m-1}\times \R$ to be the maximal open neighborhood where the map is defined and put $\nu:\mathcal N^m\to \ss^{m+1}$ to be the resulting map. 

It is easy to see that $\gamma_l(t)$ is a future directed null geodesic such that $\gamma'_l(0)$ is $g$-orthogonal to $(\mu\circ  \lambda)_*(T_{l}\mathcal L)\subset T_{\mu\circ  \lambda(l)} \ss.$ Thus $\nu$ is a mapped hypersurface corresponding to a congruence of such null geodesics. 

If $\lambda:\mathcal L^{m-1}\to \mathcal M^m$ is an immersion whose normal bundle is orientable, then there are exactly two unit lengths vector fields that are $\ov g$-orthogonal to $\lambda$ and they define two Legendrian immersions $\mathcal L\to ST\mathcal M.$ The union of the mapped hypersurfaces $\nu$ constructed for these two Legendrian immersions should be thought of as the wave front associated to $\mu\circ \lambda(\mathcal L).$

\begin{thm}\label{mappednullsurface}
Let $(\ss^{m+1},g)$ be a space-time, let $\mu:\mathcal M^m\to \ss$ be an immersed spacelike surface, and let $\wt \lambda:\mathcal L^{m-1}\to ST^*\mathcal M$ be a Legendrian map.  Let $\nu:\mathcal N^m\to \ss^{m+1}$ be the map obtained as above from $\mu$ and $\wt \lambda.$ Then the following two statements hold:
\begin{description} 
\item[1] $\nu$ is a mapped null hypersurface. In particular, the map $\wt \nu:\mathcal N\to ST\ss=ST^*\ss$ that sends $(l,t)\in \mathcal N\subset \mathcal L\times \R$ to the direction of $\gamma'_l(t)$ is a Legendrian map such that $\pr_{\ss}\circ \wt \nu=\nu,$ see Example~\ref{contactST*M}.
\item[2] If $\wt \lambda$ is a Legendrian immersion, then $\wt \nu$ also is a Legendrian immersion.
\end{description} 
\end{thm}

\pp {\it Let us prove statement $1$ of the Theorem.\/}
We have that $\mathcal L\times \{0\}\subset \mathcal N\subset \mathcal L\times \R$ and that $\mathcal N\cap (l\times \R)$ is 
connected, for all $l\in \mathcal L$. The map $\nu$ is smooth, since the velocity vectors $\gamma'_l(0)=N_l$ smoothly depend on $l\in \mathcal L$ and since $\nu(l,t)=\gamma_l(t).$

Consider the vector field $\wt N=\wt N_n=({\bf 0}, \frac{\partial}{\partial t})$ on $\mathcal N\subset \mathcal L\times \R.$ Define the vector field $N_n, n=(l,t)\in \mathcal N\subset \mathcal L\times \R,$ along $\nu$ via $N_n=\nu_*(\wt N_n)=\nu_*(l,t)({\bf 0}, \frac{\partial}{\partial t})\in T_{\nu(n)}\ss.$ Clearly $N_n=\gamma_l'(t).$ Also $N_{(l,0)}=N_l$ for all $l\in \mathcal L.$
 
Put $L_{n}\subset T_{\nu(n)}\ss$ to be the line generated by $\gamma_l'(t).$ We get the smooth line field $L_{n}, n\in \mathcal N,$ along $\nu.$
Since $\gamma_l$ are  null geodesics, all the nonzero vectors in the lines $L_{n}$ are null. Also $L_{n}\subset \nu_*(T_{n}\mathcal N)$ by construction.

Thus to prove the Theorem it suffices to show that $g\bigl (N_n, \nu_*(\wt Z_{n})\bigr)=0$ 
for all $n\in \mathcal  N, \wt Z_{n}\in T_{n}\mathcal N.$ 
Fix $n_0=(l_0,t_0)$, and $\wt Z_{n_0}\in T_{n_0}\mathcal N.$ Extend $\wt Z_{n_0}$ to a smooth vector field $\wt Z=\wt Z_n, n\in \mathcal N,$ on $\mathcal N$ such that $[\wt N, \wt Z]$ vanishes in a neighborhood of $(l_0\times \R)\cap \mathcal N.$ 

Consider the following commutative diagram: 
\begin{equation}\label{equationpullback} 
\CD 
(T\mathcal N, \wt g)@>j>> (\nu^*T\ss, \wh g, \nabla^{\wh g}) @>F>> (T\ss, g, 
\nabla^g)\\ 
@VV\tau_{\mathcal N} V @VV \nu^*\tau_{\ss}V @VV \tau_{\ss}V\\ 
\mathcal N @>\id >> \mathcal N @>\nu>> \ss . \\ 
\endCD 
\end{equation}
Here $\tau_{\ss}:T\ss\to \ss$ is the tangent bundle, $\nu^*\tau_{\ss}:\nu^*T\ss\to \mathcal N$ is the induced bundle, $\wh g=\nu^*g$ is the induced tensor field on $\nu^*T\ss$, $\nabla^{\wh g}$ is the connection on $\nu^*\tau_{\ss}$ 
induced from $\nabla^g,$ $j:T\mathcal N\to \nu^*T\ss$ is the natural bundle map, and $\wt g=\nu^* 
g$ is the induced tensor field on $T\mathcal N.$ Put $\wh N=j(\wt N)$ and $\wh Z=j(\wt Z)$ to 
be the sections of the vector bundle $\nu^*\tau_{\ss}.$ 

Let $T$ be the torsion tensor field of $\nabla^g$ and let $\wt W_i:\mathcal N\to T\mathcal N, i=1,2,$ be vector fields.  We have 
$T(\nu_*(\wt W_1), \nu_*(\wt W_2))=
\nabla^{\wh g}_{\wt W_1}j(\wt W_2)-\nabla^{\wh g}_{\wt W_2}j(\wt W_1)-j([\wt W_1, \wt W_2])$, 
see~\cite[Lemma in Section 
2.5]{GromollKlingenbergMeyer}. Since $\nabla^g$ is torsion free, we have 
\begin{equation}\label{strangetorsionfree} 
\nabla^{\wh g}_{\wt W_1}j(\wt W_2)-\nabla^{\wh g}_{\wt W_2}j(\wt W_1)=j([\wt W_1, \wt W_2]), 
\end{equation} 
for every two smooth vector fields $\wt W_1, \wt W_2:\mathcal N\to T\mathcal N.$

Since $\nabla^g$ is compatible with $g,$ we have 
\begin{equation}\label{strangemetriccompatible} 
\wt W\wh g(\wh Z_1, \wh Z_2)=\wh g \bigl(\nabla^{\wh g}_{\wt W} \wh Z_1, \wh 
Z_2\bigr)+ 
\wh g \bigl (\wh Z_1, \nabla^{\wh g}_{\wt W} \wh Z_2\bigr), 
\end{equation} 
for every vector field $\wt W:\mathcal N\to T\mathcal N$ and 
every two sections $\wh Z_1, \wh Z_2$ of $\nu^*\tau_{\ss}:\nu^*T\ss\to \mathcal N.$ 
This identity~\eqref{strangemetriccompatible} is proved in~\cite[Lemma in Section 3.4]{GromollKlingenbergMeyer} for connections induced from connections compatible with a Riemannian metric. However the same proof works for connections compatible with a Lorentzian metric.

Clearly, $g\bigl (N_{n_0}, \nu_*(\wt Z_{n_0})\bigr)=\wt g\bigl(\wt  N_{n_0}, \wt Z_{n_0}\bigr).$ 
Using identity~\eqref{strangemetriccompatible} and the fact that the vectors $N_{(l,t)}$ are the velocity vectors $\gamma'_l(t)$ of the null geodesics, we have
\begin{equation}\label{chunk1}
\wt N\wt g (\wt N, \wt Z)=\wt N\wh g(\wh N, \wh Z)=\wh g(\nabla^{\wh g}_{\wt N}\wh N, \wh Z)+ 
\wh g(\wh N, \nabla^{\wh g}_{\wt  N}\wh Z)=0+\wh g(\wh N, \nabla^{\wh g}_{\wt 
N}\wh Z).
\end{equation}

Using identities~\eqref{strangetorsionfree},~\eqref{strangemetriccompatible} and the fact that the vectors $N_n$ are null, we have 
\begin{equation}\label{chunk2}
\begin{split}
\wh g(\wh N, \nabla^{\wh g}_{\wt N}\wh Z)=\wh 
g\Bigl (\wh N, \nabla^{\wh g}_{\wt Z}\wh N+j([\wt N, \wt Z])\Bigr)=\wh g\bigl 
(\wh N, \nabla^{\wh g}_{\wt Z}\wh N+0)\bigr)=\\ 
\frac{1}{2}\wh g (\nabla^{\wh g}_{\wt Z}\wh N, \wh N)+\frac{1}{2}\wh g (\wh N, 
\nabla^{\wh g}_{\wt Z}\wh N)=\frac{1}{2}\wt Z\wh g(\wh N, \wh N)=\frac{1}{2}\wt 
Z 0=0.
\end{split}
\end{equation}

Combining equations~\eqref{chunk1} and~\eqref{chunk2} we have $\wt N\wt g(\wt N, \wt Z)=0.$ Since $\wt N=({\bf 0}, \frac{\partial}{\partial t}),$ 
we have  
\begin{equation}\label{shift}
\wt g(\wt N_{n_0}, \wt Z_{n_0})= \wt g(\wt N_{(l_0, t_0)}, \wt Z_{(l_0, t_0)})=
\wt g(\wt N_{(l_0,0)}, \wt Z_{(l_0,0)}).
\end{equation}
Decompose $\wt Z_{(l_0,0)}$ as $\wt Z^{\mathcal L}_{(l_0, 0)}+r\wt N_{(l_0, 0)},$ with 
$\wt Z^{\mathcal L}_{(l_0, 0)}\in T_{(l_0, 0)}(L\times 0).$ We identify $T_{(l_0,0)}(\mathcal L\times 0)$ 
with 
$T_{l_0}\mathcal L$ and we denote by $\wt Z_{l_0}^{\mathcal L}\in T_{l_0}\mathcal L$ the vector corresponding to $\wt 
Z^{\mathcal L}_{(l_0,0)}\in T_{(l_0,0)}(\mathcal L\times 0).$ We have 
\begin{equation}\label{step1} 
\begin{split} 
\wt g( \wt N_{(l_0,0)}, \wt Z_{(l_0,0)})=\wt g\bigl( \wt N_{(l_0,0)}, \wt 
Z^{\mathcal L}_{(l_0,0)}+r \wt N_{(l_0,0)}\bigr)= \\ r g\bigl (\nu_*(\wt N_{(l_0,0)}), \nu_*(\wt 
N_{(l_0,0)})\bigr)+ g\bigl (\nu_*(\wt N_{(l_0,0)}), \nu_*(\wt Z^{\mathcal L}_{(l_0,0)})\bigr)= \\ rg(N_{l_0}, N_{l_0})+g\bigl (N_{l_0}, (\mu\circ  \lambda)_*(\wt Z^{\mathcal L}_{l_0})\bigr)=0+g\bigl (N_{l_0}, (\mu\circ  \lambda)_*(\wt Z^{\mathcal L}_{l_0})\bigr). 
\end{split} 
\end{equation} 

Recall that $N_{l_0}=(N_{l_0}^s, N_{l_0}^t)\in \mu_*(T_{ \lambda(l_0)}\mathcal M)\oplus (\mu_*(T_{ \lambda(l_0)}\mathcal M))^{\perp}=T_{\mu\circ  \lambda(l_0)}\ss$ and that $N_{l_0}^s=\mu_*(X_{l_0}),$ where $X_{l_0}\in T_{ \lambda(l_0)}M$ is the unit vector whose equivalence class is $\wt \lambda(l_0).$
Thus 
\begin{equation}\label{step2}
\begin{split}
g\bigl (N_{l_0}, (\mu\circ  \lambda)_*(\wt Z^{\mathcal L}_{l_0})\bigr)=g\Bigl 
(N^s_{l_0}+ N^t_{l_0}, \mu_*( \lambda_*(\wt Z^{\mathcal L}_{l_0}))\Bigr)=\\
g\bigl 
(\mu_*(X_{l_0}), \mu_*( \lambda_*(\wt Z^{\mathcal L}_{l_0}))\bigr)+g\bigl 
(N^t_{l_0}, \mu_*( \lambda_*(\wt Z^{\mathcal L}_{l_0}))\bigr)=\ov g\bigl (X_{l_0},  \lambda_*(\wt Z^{\mathcal L}_{l_0})\bigr)+0.
\end{split}
\end{equation}
Since $\wt \lambda$ is Legendrian, $X_{l_0}$ is $\ov g$-orthogonal to $\lambda_*(T_{l_0}\mathcal L)\subset T_{ \lambda(l_0)}\mathcal M$ and hence $\ov g\bigl (X_{l_0},  \lambda_*(\wt Z^{\mathcal L}_{l_0})\bigr)=0.$
Combining equations~\eqref{shift},~\eqref{step1},~\eqref{step2} we have 
\begin{equation}
\begin{split}
g(N_{n_0}, Z_{n_0})=\wt g(\wt N_{(l_0, t_0)}, \wt Z_{(l_0, t_0)})=\wt g(\wt N_{(l_0, 0)}, \wt Z_{(l_0, 0)})=\\
g\bigl (N_{l_0}, (\mu\circ  \lambda)_*(\wt Z^{\mathcal L}_{l_0})\bigr)=\ov g\bigl (X_{l_0},  \lambda_*(\wt Z^{\mathcal L}_{l_0})\bigr)=0.
\end{split}
\end{equation}
This finishes the proof of Statement $1$ of the Theorem.

{\it Let us prove statement $2.$\/} Here the main difficulty is that even when $(\ss, g)$ is geodesically complete, the geodesic flow on $T\ss$ does not seem to give rise to a flow on $ST\ss$ or on the subspace of it formed by the null directions, except for some very special $(\ss, g).$

Consider the map $\wt{\exp'}:T\ss\to T\ss$ that associates to $v\in T_{x}\ss$ the velocity vector $\gamma'_v(1)\in T_{\gamma_v(1)}\ss$ of the unique inextendible geodesic $\gamma_v(t)$ with $\gamma_v(0)=x$ and $\gamma'_v(0)=v.$ Put $U'\subset T\ss$ to be the (maximal) domain of this map. It is an open set, see~\cite[discussion after Lemma 1 in Section 2.8 and Proposition in Section 2.9]{GromollKlingenbergMeyer}. Clearly $\wt{\exp'}:U'\to U'$ is a smooth bijection. The inverse map sends $v\in T_x\ss$ to $\wt {\exp'}(-v)$ and hence is also smooth. Thus $\wt{\exp'}:U'\to U'$ is a diffeomorphism.

Put $O\subset U'\subset T\ss$ to be (the image of) the zero section of $T\ss\to \ss.$
Put $U=U'\setminus O.$ Clearly the restriction $\wt{\exp'}|_U$ is a diffeomorphism $U\to U$ that we denote by $\wt \exp.$

Consider the map $\wt \kappa:\mathcal L\to ST\ss$ that is described by the pair: the map $\kappa=\mu\circ \lambda:\mathcal L\to \ss$ and the vector field $N_l\in T_{\kappa(l)}\ss, l\in \mathcal L,$ along $\kappa.$ Let us show that $\wt \kappa$ is an immersion. Take $l\in \mathcal L$ and its neighborhood $\mathcal O\subset \mathcal L$ such that $\lambda(\mathcal O)$ is contained in an open neighborhood $\mathcal P\subset \mathcal M$ for which the restriction $\mu|_{\mathcal P}:\mathcal P\to \ss$ is an embedding. It suffices to show that $\wt \kappa|_{\mathcal O}$ is an immersion. The restriction of the bundle $\pr_{\mathcal M}:ST\mathcal M\to \mathcal M$ to $\mathcal P\subset \mathcal M$ gives the $S^{m-1}$-bundle $ST\mathcal P\to \mathcal P.$ The restriction of $\pr_{\ss}:ST\ss\to \ss$ to $\mu(\mathcal P)$ gives the $S^m$-bundle $ST\ss|_{\mu(\mathcal P)}\to \mu(\mathcal P).$ 
The embedding $\mu|_{\mathcal P}$ induces the natural bundle map 
\begin{equation}\label{equation1} 
\CD 
ST\mathcal P@>i>> ST\ss|_{\mu(\mathcal P)}\\ 
@VVV @VVV\\ 
\mathcal P @>\mu|_{\mathcal P} >> \mu(\mathcal P). \\ 
\endCD 
\end{equation}
For $p\in \mathcal P$ put $V_{p}\in T_{\mu(p)}\ss$ to be the unique future pointing timelike vector such that $g(V_p, V_p)=-1$ and $V_p$ is $g$-orthogonal to $\mu_*(T_p\mathcal M).$ Put $V, -V\subset ST\ss|_{\mu(\mathcal P)}$ to be the images of the two sections of $ST\ss|_{\mu(\mathcal P)}\to \mu(\mathcal P)$ that send $\mu(p)$ to the direction of $V_p$ and to the direction of $-V_p,$ respectively. The direct sum decomposition $T_{\mu(p)}\ss=\mu_*(T_p\mathcal P)\oplus \span(V_p)$ induces the natural fiber preserving smooth map $\pi:ST\ss|_{\mu(\mathcal P)}\setminus (V\sqcup -V)\to ST\mathcal P.$
For all $l\in \mathcal O$ we have $\mu_*(X_l)+V_{\mu(l)}=N_l.$ Thus the maps $i\circ \wt \lambda|_{\mathcal O}$ and $\pi \circ \wt \kappa|_{\mathcal O}:\mathcal O\to ST\mathcal P$ are equal. Since $i$ is an embedding, $\wt \lambda$ is an immersion, and $\pi$ is smooth, we get that $\wt \kappa|_{\mathcal O}$ is an immersion, and hence $\wt \kappa$ is an immersion. Put $\wh \kappa:\mathcal L\to T\ss\setminus O$ to be the map described by the pair $\kappa$ and the vector field $N_l, l\in \mathcal L,$ along $\kappa.$ Since $\wt \kappa$ is immersion and it is a composition of $\wh \kappa$ and the smooth quotient map $T\ss\setminus O\to ST\ss,$ we get that $\wh \kappa$ is an immersion.

Define the map $\wh \nu:\mathcal N\to T\ss$ by sending $(l, t)\in \mathcal N\subset \mathcal L\times \R$ to $\gamma'_{l}(t)\in T_{\nu(l,t)}\ss.$ Let us show that $\wh \nu$ is an immersion. Put $\mathcal N^+, \mathcal N^-, \mathcal N^0\subset \mathcal N\subset \mathcal L\times \R$ to be the subsets formed by points $(l,t)$ whose $t$ coordinate is respectively greater than zero, less than zero, and is equal to zero.

Clearly $\mathcal N^+$ is open and $\wh \nu(l,t)=\wt \exp(t N_l),$
for all $(l,t)\in \mathcal N^+.$ Since $\wh \kappa:\mathcal L\to T\ss\setminus O$ is an immersion, we get that the map $\beta^+:\mathcal N^+\to T\ss\setminus O$ that sends $(l,t)\in \mathcal N^+$ to $tN_l$ is an immersion. Since $\wt\exp$ is a diffeomorphism, we get that $\wh\nu$ is an immersion at all points of $\mathcal N^+.$ Similarly one gets that $\wh \nu$ is an immersion at all points of $\mathcal N^-.$ 

Take $(l,0)\in N^0$ and a nonzero tangent vector $(v_{\mathcal L}, v_{\R})\in T_{(l,0)}\mathcal N\subset T_{(l,0)}(\mathcal L\times \R)=T_{l}\mathcal L\oplus T_0\R=T_{l}\mathcal L\oplus \R.$ 
Then $\wh\nu_*(v_{\mathcal L}, v_{\R})\in T_{N_l} (T\ss)$ and 
$(\tau_{T\ss})_*\circ \wh\nu_*(v_{\mathcal L}, v_{\R})=(\tau_{T\ss})_*\circ \wh \kappa_*(v_{\mathcal L})+v_{\R}N_l\in T_{\nu(l,0)}\ss.$ Since $(\tau_{T\ss})_*\circ \wh \kappa_*(v_{\mathcal L})\in \mu_*(T_{\lambda(l)}\mathcal M)$ and $N_l\not \in \mu_*(T_{\lambda(l)}\mathcal M),$
we get that $\wh \nu_*(v_{\mathcal L}, v_{\R})\neq 0$ if $v_{\R}\neq 0.$ On the other hand $\wh \nu_*(v_{\mathcal L}, 0)=\wh \kappa(v_{\mathcal L})$ is nonzero since $\wh \kappa$ is an immersion. Thus $\wh \nu$ is an immersion at all the points of $\mathcal N^0.$

Let $q:T\ss\setminus O\to ST\ss$ be the quotient map by the action of $\R^+$ that we used to define $ST\ss.$ Clearly $\wt \nu=q\circ \wh \nu.$ Since $\wh \nu$ is an immersion, to prove that $\wt \nu$ is an immersion it suffices to show that for every $\ov n\in \mathcal N$ and nonzero $v\in T_{\ov n}\mathcal N$ the nonzero vector $(\wh \nu)_*(v)$ is not tangent to the $\R^+$-fiber of $q$ containing $\wh \nu(\ov n).$

We prove this by considering three cases: $\ov n\in \mathcal N^+, \ov n\in \mathcal N^-,$ and $\ov n\in \mathcal N^0.$

Assume that $\ov n=(\ov l,\ov t)\in \mathcal N^+$ and that $\wh \nu_*(v)$ is tangent to the $\R^+$-fiber of $q$ containing $\wh \nu(\ov n).$ Let $\alpha:(-\epsilon, \epsilon)\to T\ss\setminus O$ defined by $\alpha(\tau)=\wh \nu(\ov n)+\tau \wh \nu(\ov n)$ be the parameterization of a small part of the $\R^+$-fiber of $q$ that contains $\wh \nu(\ov n).$ Since $\wt{\exp}$ is a diffeomorphism, we get that $(\wt {\exp})^{-1}_*\circ \wh \nu_*(v)$ is a nonzero vector tangent to the curve $\wt \alpha=\wt {\exp}^{-1}\circ \alpha$ at  $\wt {\exp}^{-1}\circ \alpha(0).$

Let $\gamma(t)$ be the null geodesic such that $\gamma(0)=\mu\circ \lambda(\ov l)$ and $\gamma'(0)=\ov tN_{\ov l}.$ Since $\wh \nu|_{\mathcal N^+}=\wt \exp\circ \beta^+$ we get that $\gamma(1)=\tau_{\ss}(\wh \nu(\ov n))$ and $\gamma'(1)=\wh \nu(\ov n).$ From the definition of $\wt \exp$ we get that $\wt \alpha(\tau)=\wt \exp^{-1}(\alpha(\tau))=(\tau+1)\gamma'(-\tau)\in T_{\gamma(-\tau)}\ss,$ for all $\tau\in (-\epsilon, \epsilon).$

Now $\wt {\exp}^{-1}_*\circ \wh \nu_*(v)$ is a nonzero vector tangent to the immersed submanifold $S=\{tN_l|t\in \R, l\in \mathcal L\}\subset T\ss$ at the point $\ov tN_{\ov l}.$ Since $\tau_{\ss}(tN_l)=\mu\circ \lambda(l)$ for all $t\in \R, l\in \mathcal L,$ 
we get that $(\tau_{\ss})_*\bigl(\wt {\exp}^{-1}_*\circ \wh \nu_*(v)\bigr)\in \mu_*(T_{\lambda(\ov l)}\mathcal M).$
Clearly $(\tau_{\ss})_*(\wt \alpha'(0))=-\gamma'(0)=-\ov t N_{\ov l}\not \in \mu_*(T_{\lambda(\ov l)}\mathcal M).$ Thus $\wt {\exp}^{-1}_*\circ \wh \nu_*(v)$ is not tangent to $\wt \alpha$ at $\wt \alpha(0)$ and $\wt \nu$ is an immersion at $\ov n\in \mathcal N^+.$

Hence $\wt \nu$ is an immersion at all the points of $\mathcal N^+.$
Similarly one gets that $\wt \nu$ is an immersion at all the points of $\mathcal N^-.$

Let $\ov n=(\ov l, 0)$ be a point of $\mathcal N^0$ and let $(v_{\mathcal L}, v_{\R})\in T_{(\ov l, 0)}\mathcal N=T_{\ov l}\mathcal L\oplus T_{0}\R=T_{\ov l}\mathcal L\oplus \R$ be a nonzero tangent vector. Let us show that $\wh \nu_*(v_{\mathcal L}, v_{\R})$ is not tangent to the $\R^+$-fiber of $q$ containing $\wh \nu(\ov l, 0)=N_{\ov l}.$ Note that $(\tau_{T\ss})_*$ applied to any vector tangent to the $\R^+$-fiber of $q$ is zero, while, as we discussed above, $(\tau_{T\ss})_*\circ \wh \nu_*(v_{\mathcal L}, v_{\R})=(\tau_{T\ss})_*\circ \wh \kappa_*(v_{\mathcal L})+v_{\R}N_{\ov l}$ is nonzero for $v_{\R}\neq 0.$ This give the proof for vectors $(v_{\mathcal L}, v_{\R})$ with nonzero $v_{\R}.$

Note that $\wh \nu_*(v_{\mathcal L}, 0)=\wh \kappa_*(v_{\mathcal L}).$ Clearly $\wt \kappa=q\circ \wh \kappa$ and since $\wt \kappa$ is an immersion we get that $q_*\circ \wh \kappa_*(v_{\mathcal L})=q_*\circ \wh \nu_*(v_{\mathcal L}, 0)$ is nonzero for every nonzero $v_{\mathcal L}.$ On the other hand, $q_*$ applied to any vector tangent to the $\R^+$-fiber of $q$ is zero. Thus $\wt \nu$ is an immersion at all the points of $\mathcal N^0\subset \mathcal N$ and hence $\wt \nu:\mathcal N\to ST\ss$ is an immersion.
\qed

\begin{rem}\label{trivial}
Let $\nu:\mathcal N^m\to \ss^{m+1}$ be a mapped null hypersurface and let $h:\mathcal N'\to \mathcal N$ be a diffeomorphism. Then clearly $\nu\circ h:\mathcal N'\to \ss$ is a mapped null hypersurface. 

If the natural map $\wt \nu:\mathcal N^m\to ST\ss=ST^*\ss$ is an immersion, then it is a Legendrian immersion and the map $\wt {\nu \circ h}:\mathcal N'\to ST\ss=ST^*\ss$ associated with the mapped null hypersurface $\nu \circ h:\mathcal N'\to \ss$ also is a Legendrian immersion.

Similarly if $U\subset \mathcal N$ is open, then $\nu|_U:U\to \ss$ is a mapped null hypersurface. Note that if $\nu|_U$ is an embedding, then $\nu|_U:U\to \ss$ is an embedded null hypersurface.
\end{rem}

\section{From mapped null hypersurfaces to Legendrian maps.}
Let $(\ss^{m+1}, g)$ be a space-time. Let $\nu:\mathcal N^m\to \ss^{m+1}$ be a mapped null hypersurface, let $\mu:\mathcal M^m\to \ss^{m+1}$ be an immersed spacelike hypersurface, and let  $\mathcal L_{\mu, \nu}$ be the pull back of the maps $\mu$ and $\nu.$ We will show that $\mu$ and $\nu$ 
canonically define a Legendrian map $\wt \lambda_{\mu, \nu}: \mathcal L_{\mu, \nu}\to ST^*\mathcal M$ of the $(m-1)$-dimensional pull-back manifold and that $\Im (\pr_{\mathcal M}\circ \wt \lambda_{\mu, \nu})=\Im \mu^{-1}(\Im \mu\cap \Im \nu).$ 

We will also show that if the map $\wt \nu:\mathcal N^m\to ST\ss$ associated to $\nu$ is an immersion, then $\wt \lambda_{\mu, \nu}$ is a Legendrian immersion. In this case the singularities of $\pr_{\mathcal M}\circ \wt \lambda_{\mu, \nu}$ are Legendrian singularities. In particular, this is so when $\nu$ is the mapped null hypersurface arising from a Legendrian immersion $\wt \lambda':\mathcal L'\to ST^*\mathcal M'$ for some immersed spacelike hypersurface $\mu':\mathcal M'\to \ss,$ see Theorem~\ref{mappednullsurface}.

The Lorentz metric $g$ allows us to identify $ST\ss$ with $ST^*\ss.$ 
Let $N_n\in T_{\nu(n)}\ss, n\in \mathcal N,$ 
be a smooth nowhere zero null vector field 
along $\nu$ such that for all $n\in \mathcal N$ the equivalence class of $N_n$ is $\wt \nu(n)\in ST\ss=ST^*\ss.$ For $n\in \mathcal N,$ put $\theta_n\in T^*_{\nu(n)}\ss$ to be the nonzero covector such that $\theta_n(v)=g(N_n, v),$ for all $v\in T_{\nu(n)}\ss.$ We get the smooth nowhere zero covector field $\theta_n, n\in \mathcal N,$ along $\nu$ such that for all $n\in \mathcal N$ the equivalence class of $\theta_n$ is $\wt \nu(n)\in ST^*\ss=ST\ss.$

Consider the pull-back diagram 
\begin{equation}\label{pullbacksquare} 
\CD 
\mathcal L_{\mu, \nu}@> \lambda_{\mu, \nu} >> \mathcal M ^m\\ 
@VV j V @VV \mu V\\ 
\mathcal N^{m} @> \nu >> \ss^{m+1}.\\ 
\endCD 
\end{equation} 

By definition of the pull-back $\mathcal L_{\mu, \nu}=\{(m, n)\in \mathcal M\times \mathcal N|\mu(m)=\nu(n)\}\subset \mathcal M\times \mathcal N.$
Choose $(m,n)\in \mathcal L_{\mu, \nu}.$  Since $\mu$ is an immersion, $\mu_*(T_m\mathcal M)$ is $m$-dimensional. Since $\nu$ is a mapped null hypersurface and by definition of $\wt \nu,$ the nonzero vector $N_n\in \nu_*(T_n\mathcal N)$ is null. Since $\mu$ is spacelike, all the nonzero vectors in $\mu_*(T_m\mathcal M)$ are spacelike, and hence $N_n\not \in \mu_*(T_m\mathcal M).$ For dimension reasons we get that the minimal linear subspace of $T_{\nu(n)}\ss=T_{\mu(m)}\ss$ that contains $\mu_*(T_m\mathcal M)\cup\nu_*(T_n\mathcal N)$ is equal to $T_{\nu(n)}\ss=T_{\mu(m)}\ss.$ Thus $\mu$ and $\nu$ are transverse and hence $\mathcal L_{\mu, \nu}$ is an $(m-1)$-dimensional smooth embedded submanifold of $\mathcal M\times \mathcal N.$

Clearly $\lambda_{\mu, \nu}(m,n)=m$ and $j(m,n)=n,$ for $(m,n)\in \mathcal L_{\mu, \nu}.$ We define the smooth covector field $\phi_l\in T^*_{\lambda(l)}\mathcal M, l\in \mathcal L_{\mu, \nu},$ along $\lambda_{\mu, \nu}$ as follows. For $l=(m,n)\in \mathcal L_{\mu, \nu}$ and $v\in T_{\lambda(l)}\mathcal M=T_m\mathcal M$ put $\phi(v)=\theta_n(\mu_*(v)).$ Recall that the covector $\theta_n$ is nonzero, $\theta_n|_{\nu_*(T_n\mathcal N)}$ is zero, and $T_{\mu(m)}\ss=T_{\nu(n)}\ss$ is the linear span of $\mu_*(T_m\mathcal M)\cup \nu_*(T_n\mathcal N).$ Thus the covector field $\phi_l, l\in \mathcal L_{\mu, \nu},$ along $\lambda_{\mu, \nu}$ is nowhere zero. Hence the pair: 
$\lambda_{\mu, \nu}$ and the covector field $\phi_l, l\in \mathcal L_{\mu, \nu},$ along $\lambda_{\mu, \nu}$ define a map $\wt \lambda_{\mu, \nu}:\mathcal L_{\mu, \nu}\to ST^*\mathcal M=ST\mathcal M.$ It is easy to see that the map $\wt \lambda_{\mu, \nu}$ does not depend on the choice of the vector field $N_n, n\in \mathcal N,$ along $\nu$ from which we started the construction.

\begin{thm}\label{tripback}
Let $(\ss^{m+1}, g)$ be a space-time, let $\nu:\mathcal N^m\to \ss^{m+1}$ be a mapped null hypersurface, and let $\mu:\mathcal M^m\to \ss^{m+1}$ be an immersed spacelike hypersurface. Let $\mathcal L_{\mu, \nu}$ be the smooth $(m-1)$-dimensional manifold that is the pull-back of $\mu$ and $\nu.$ Let  
$\wt \lambda_{\mu, \nu}:\mathcal L_{\mu, \nu}\to ST^*\mathcal M$ be the map constructed above. 
Then 
\begin{description}
\item[1] The map $\wt \lambda_{\mu, \nu}:\mathcal L_{\mu, \nu}\to ST^*\mathcal M$ is a Legendrian map and $\Im (\pr_{\mathcal M}\circ \wt \lambda_{\mu, \nu})=\Im \mu^{-1}(\Im \mu\cap \Im \nu).$
 \item[2] If the map $\wt \nu:\mathcal N\to ST\ss$ that is naturally associated with $\nu$ is an immersion, then the map $\wt \lambda_{\mu, \nu}:\mathcal L_{\mu, \nu}\to ST^*\mathcal M$ is a Legendrian immersion. 
\end{description}
\end{thm}

\pp {\it Let us prove statement $1$ of the Theorem.\/}
The fact that $\Im (\pr_{\mathcal M}\circ \wt \lambda_{\mu, \nu})=\Im \mu^{-1}(\Im \mu\cap \Im \nu)$ is clear from the construction of $\wt \lambda_{\mu, \nu}.$ 

To see that $\wt \lambda_{\mu, \nu}$ is a Legendrian map it suffices to show that $\phi_l\bigl((\lambda_{\mu, \nu})_*(v)\bigr)=0$ for every $l\in \mathcal L_{\mu, \nu}$ and $v\in T_l\mathcal L_{\mu, \nu}.$ 

By definition of $\phi_l$ we have $\phi_l\bigl((\lambda_{\mu, \nu})_*(v)\bigr)=\theta_n\bigl((\mu\circ \lambda_{\mu, \nu})_*(v)\bigr).$ Since the diagram~\eqref{pullbacksquare} is commutative, we have $\theta_n\bigl((\mu\circ \lambda_{\mu, \nu})_*(v)\bigr)=\theta_n\bigl(\nu_*(j_*(v))\bigr).$ By definition of $\wt \nu$ we have that $\theta_n|_{\Im \nu_*(T_n\mathcal N)}=0.$ Thus $\phi_l\bigl((\lambda_{\mu, \nu})_*(v)\bigr)=\theta_n\bigl(\nu_*(j_*(v))\bigr)=0.$

{\it To prove statement $2$ of the Theorem\/} we will show that for every $\ov l=(\ov m, \ov n)\in \mathcal L_{\mu, \nu}$ the map $\wt \lambda_{\mu, \nu}$ is an immersion at $\ov l.$ Put $\mathcal P\subset \mathcal M$ to be an open neighborhood such that $\ov m\in \mathcal P$ and $\mu|_{\mathcal P}$ is an embedding. 

Put $\mathcal L=\mathcal L_{\mathcal P}=\{(m,n)\in \mathcal M\times \mathcal N|\mu(m)=\nu(n) \text{ and }m\in \mathcal P\}\subset \mathcal L_{\mu, \nu}.$ We will denote $\wt \lambda_{\mu, \nu}|_{\mathcal L}:\mathcal L\to ST^*\mathcal P\subset ST\mathcal M$ by $\wt \lambda_{\mathcal L}$ and we will denote $\lambda_{\mu, \nu}|_{\mathcal L}:\mathcal L\to \mathcal P$ by $\lambda_{\mathcal L}.$ It suffices to show that $\wt \lambda_{\mathcal L}$ is an immersion at $\ov l\in \mathcal L.$ 

Take a nonzero vector $v=(v_{\mathcal P}, v_{\mathcal N})\in T_{\ov l}\mathcal L\subset T_{(\ov m, \ov n)}(\mathcal P\times \mathcal N)=T_{\ov m}\mathcal P\oplus T_{\ov n}\mathcal N.$ From the construction of $\lambda_{\mathcal L}$ and $\wt \lambda_{\mathcal L}$ we get that $(\pr_{\mathcal P})_*\circ (\wt \lambda_{\mathcal L})_*(v)=v_{\mathcal P}.$ Thus $(\wt \lambda_{\mathcal L})_*(v)\neq {\bf 0}$ if $v_{\mathcal P}\neq {\bf 0}.$ Hence it suffices to show that $(\wt \lambda_{\mathcal L})_*({\bf 0}, v_{\mathcal N})\neq {\bf 0}$ for $({\bf 0}, v_{\mathcal N})\in T_{\ov l}\mathcal L$ with nonzero 
$v_{\mathcal N}.$

The embedding $\mu|_{\mathcal P}:\mathcal P\to \ss$ induces the diffeomorphism
$ST\mu|_{\mathcal P}:ST\mathcal P\to ST\mu(\mathcal P)$ onto the spherical tangent bundle of $\mu(\mathcal P).$ 

Consider the map $\wt \nu\circ j:\mathcal L\to ST\ss$ that maps $(m,n)\in \mathcal L$ to the equivalence class of $N_n\in T_{\nu(n)}\ss.$ Since $\mathcal L$ is the pull-back of $\nu$ and $\mu|_{\mathcal P}$ we get that $\nu(n)\in \mu(\mathcal P)$ for all $(m,n)\in \mathcal L$. Thus $\wt \nu \circ j(\mathcal L)$ is in the total space $ST\ss|_{\mu(\mathcal P)}$ of the restriction of the $S^m$-bundle $ST\ss\to \ss$ to $\mu(\mathcal P)\subset \ss.$ Consider the $S^{m-1}$-bundle $\mathcal C\to \ss$ whose total space $\mathcal C\subset ST\ss$ is formed by the future pointing null directions. Clearly 
$\wt \nu(\mathcal N)\subset \mathcal C.$ Thus $\wt \nu\circ j(\mathcal L)$ is in the total $\mathcal C|_{\mu(\mathcal P)}$ of the restriction of the bundle $\mathcal C\to \ss$ to $\mu(\mathcal P).$ 

Put $T\ss|_{\mu(\mathcal P)}$ to be the total space of the restriction to $\mu(\mathcal P)$ of the bundle $T\ss\to \ss.$
The $g$-orthogonal projection $T\ss|_{\mu(\mathcal P)}\to T\mu(\mathcal P)$ induces the diffeomorphism $\delta:\mathcal C|_{\mu(\mathcal P)}\to ST\mu(\mathcal P).$

The Riemannian metric $\ov g=\mu^*g$ on $T\mathcal P$ allows us to identify $ST \mathcal P$ with $ST^*\mathcal P.$ 
For $l=(m,n)\in \mathcal L$ put $V_l\in T_{\lambda_{\mathcal L}(l)}\mathcal P=T_{m}\mathcal P$ to be the unique vector such that $\phi_l(w)=\ov g(V_l, w),$
for all $w\in T_m\mathcal P.$ From the construction of $\phi_l$ it is easy to see that $\mu_*(V_l)$ is the $g$-orthogonal projection of $N_n\in T_{\nu(n)}\ss$ to 
$\mu_*(T_m\mathcal P)\subset T_{\nu(n)}\ss.$ One verifies that the equivalence class of $V_l$ in $ST\mathcal P=ST^*\mathcal P$ is $\wt \lambda_{\mathcal L}(l).$

Thus we have that  the maps $\wt\lambda_{\mathcal L}:\mathcal L\to ST\mathcal P=\mathcal ST^*P$ and $(ST\mu|_{\mathcal P})^{-1}\circ \delta\circ \wt \nu\circ j|_{\mathcal L}:\mathcal L\to ST\mathcal P$ are equal. Thus if $v=({\bf 0}, v_{\mathcal N})\in T_{\ov l}\mathcal L$ is a nonzero vector, then we have
$(\wt \lambda_{\mathcal L})_*(v)= (ST\mu|_{\mathcal P})^{-1}_*\circ \delta_*\circ \wt \nu_*\circ (j|_{\mathcal L})_*({\bf 0}, v_{\mathcal N})=(ST\mu|_{\mathcal P})^{-1}_*\circ \delta_*\circ \wt \nu_*(v_{\mathcal N}).$ Since $v_{\mathcal N}\neq {\bf 0},$ $\wt \nu$ is an immersion, and $(ST\mu_{\mathcal P})^{-1}, \delta$ are diffeomorphisms, we get that $(\wt \lambda_{\mathcal L})_*(v)\neq {\bf 0}.$
Hence $\wt \lambda_{\mu, \nu}$ is an immersion.
\qed

\begin{ex}[Null cone]
Let $(\ss, g)$ be a space-time. For $x\in \ss$ put $C^+_x$ (respectively $C^-_x$) to be the hemicone of future pointing (respectively past pointing) null vectors in $T_x\ss.$ Put $\mathcal C^+_x\subset C^+_x$ and $\mathcal C^-_x\subset C^-_x$ to be the maximal open subsets on which $\exp_x$ is well defined.

Choose a (possibly small) immersed spacelike hypersurface $\mu:\mathcal M^m\to \ss^{m+1}$ such that 
$x=\mu(\ov x)$ for some $\ov x\in \mathcal M$ and let $\ov g$ be the induced Riemannian metric on $\mathcal M.$
Let $\wt \lambda:S^{m-1}\to ST^*\mathcal M$ be a Legendrian embedding whose image is the $S^{m-1}$-fiber over the point $\ov x.$ 
Let $\nu:\mathcal N\to \ss$ be the mapped null hypersurface from Theorem~\ref{mappednullsurface} constructed using the above $\mu,\wt \lambda$ and $\mathcal L=S^{m-1}.$ By Theorem~\ref{mappednullsurface} the natural map $\wt\nu:\mathcal N\to ST\ss$ is a Legendrian immersion.

\def\mfootu{\footnote{In the work of Lerner~\cite[Lemma 2]{Lerner} it is proved that the exponential of the future null cone of $x$ is an embedded null hypersurface when restricted to the preimage under $\exp_x$ of a simple neighborhood of $x.$ We did not find more general statements about null cones giving rise to embedded null hypersurfaces in the literature. Miguel~Sanchez pointed to us that the fact that $\exp_x|_U$ is an embedded null hypersurface  also follows from the Gauss Lemma for Lorentzian manifolds~\cite{ONeill} and we thank him for this remark.}}

Let $N_l\in T_{\mu\circ  \lambda(l)}\ss=T_x\ss, l\in S^{m-1},$ be the future pointing null vector field  along $\mu\circ  \lambda$ that we used to construct $\nu.$
Consider the map $h:S^{m-1}\times \R\to (C^+_x\sqcup C^-_x\sqcup {\bf 0})\subset T_x\ss$ defined by $h(l,t)=tN_l.$ Put $\mathcal N^+$ (respectively $\mathcal N^-$) to be the open subset of $\mathcal N\subset S^{m-1}\times \R$ consisting of all the points with the positive (respectively negative) $\R$-coordinate. Clearly $h:\mathcal N^+\to \mathcal C^+$ and $h:\mathcal N^-\to \mathcal C^-$ are diffeomorphisms and 
$\exp_x(h(l,t))=\nu(l,t),$ for all $(l,t)\in \mathcal N^{\pm}.$

Combining this with Remark~\ref{trivial} we get that
$\exp_x:\mathcal  C^+_x\to \ss$ and $\exp_x:\mathcal C^-_x\to \ss$ are mapped null hypersurfaces, i.e. the exponential of the future and of the past null hemicones at $x$ are mapped null hypersurfaces on the maximal open subsets where they are defined.
In particular if an open $U\subset \mathcal C^{\pm}$ is such that $\exp_x|_U$ is an embedding, then $\exp_x:U\to \ss$ is  an embedded null hypersurface.\protect\mfootu

Moreover by Theorem~\ref{mappednullsurface} the natural maps $\mathcal C^+_x\to ST\ss=ST^*\ss$ and $\mathcal C^-_x\to ST\ss=ST^*\ss$ are Legendrian immersions.

Let $\mu':\mathcal M'\to \ss$ be an immersed spacelike hypersurface. Then by Theorem~\ref{tripback} the map $\exp_x|_{\mathcal C^+_x}:\mathcal C^+_x\to \ss$
defines the Legendrian immersion $\wt \lambda_{\mu', \exp_x|_{\mathcal C^+_x}}:\mathcal L_{\mu', \exp_x|_{\mathcal C^+_x}}\to ST\mathcal M'$ such that $\Im (\pr_{\mathcal M'}\circ \wt \lambda_{\mu', \exp_x|_{\mathcal C^+_x}})=\exp_x(\mathcal C^+_x)\cap \mu'(\mathcal M').$ Thus the intersection of the future null cone of $x$ with the spacelike immersed hypersurface $\mu(\mathcal M')$ is naturally parameterized by the projection to $\mathcal M'$ of the Legendrian immersion to $ST^*\mathcal M'.$ 

Similarly one get that the intersection of the past null cone of $x$ with the spacelike immersed hypersurface $\mu(\mathcal M')$ is also naturally parameterized by the projection to $\mathcal M'$ of a Legendrian immersion to $ST^*\mathcal M'.$ 
\end{ex}

\appendix

\section{Low's~\cite{LowLegendrian} results on null congruences and Legendrian submanifolds of the space of null geodesics in globally hyperbolic and strongly causal $(\ss^{3+1}, g).$} 

Recall a few more Lorentzian geometry definitions and facts. 
An open set in $(\ss^{m+1}, g)$ is {\em causally convex\/} if its intersection with every non-spacelike curve is connected and  $(\ss, g)$ is {\em 
strongly causal\/} if every point in it has arbitrarily small causally convex 
neighborhoods. A strongly causal space-time $(\ss, g)$ is {\it globally hyperbolic\/} if for every $x_1, x_2\in \ss$ the set of all $x\in \ss$ such that there exists a piecewise smooth nonspacelike curve from $x_1$ to $x_2$ through $x$ is compact. 

A {\em Cauchy surface\/} $M$ is a subset of a space-time $\ss$ such that for 
every inextendible non-spacelike curve $\gamma(t)$ in $\ss$ there exists 
exactly one $t_0\in \R$ with $\gamma(t_0)\in M$.  
A space-time is globally hyperbolic if and only if it admits a 
Cauchy surface, see~\cite[pages 211-212]{HawkingEllis}. 
Geroch~\cite{Geroch} showed that globally hyperbolic $(\ss, g)$ are rather simple topologically and they are homeomorphic to a product of $\R$ and a Cauchy surface. Bernal and Sanchez~\cite{BernalSanchez} showed that every globally hyperbolic space-time $(\ss^{m+1}, g)$ admits a smooth spacelike Cauchy surface and moreover $\ss$ is in fact diffeomorphic to a product
of $\R$ and this Cauchy surface.

Put $\frak N=\frak N_{(\ss, g)}$ to be the space of all null geodesics in $(\ss, g)$ up to an affine reparameterization. In general $\frak N$ is not a manifold. However 
for globally hyperbolic $(\ss, g),$ the space $\frak N$ is a smooth contact manifold 
contactomorphic to the spherical cotangent bundle $ST^*M$ of a smooth spacelike Cauchy surface $M^m\subset \ss^{m+1}.$ This fact was proved by Low~\cite[Corollary 1,
Lemma 2, Corollary 2]{LowLegendrian} for $(3+1)$-dimensional globally-hyperbolic $(\ss^{3+1}, g)$. This result and the techniques, Low used to get it, generalize to globally hyperbolic space-times of all dimensions, see Natario and Tod~\cite[pages 252-253]{NatarioTod}.

Since the Cauchy surface $M$ is spacelike we can identify $STM$ and $ST^*M.$ Under the contactomorphism $\frak N\to ST^*M$ a null geodesic $\gamma$ is mapped to the point of $ST^*M=STM$ that is the direction of the $g$-orthogonal projection to $M$ of the velocity vector of $\gamma$ at the intersection point of $\gamma$ with $M.$

Low~\cite{LowLegendrian} observed strong and fascinating relations between null congruences and Legendrian submanifolds of $\frak N$ for $(3+1)$-dimensional globally hyperbolic $(\ss^{3+1}, g).$ 
The combination of his~\cite[Lemma 2, Corollary 3]{LowLegendrian} says that the null congruences orthogonal to a $2$-dimensional spacelike surface are exactly the Legendrian submanifolds of $ST^*M=\frak N.$ Unfortunately, if taken literally this statement is false for rather technical reasons. 

For example, in order for the Legendrian submanifold to be embedded, rather than immersed, one has to require that no two points of the $2$-dimensional spacelike surface $\Sigma$ belong to the same null geodesic that is $g$-orthogonal to $\Sigma.$ This would follow automatically if the $2$-dimensional spacelike surface $\Sigma$ is a subset of some Cauchy surface. However it is easy to construct examples of $2$-dimensional embedded spacelike surfaces $\Sigma$ in globally hyperbolic $(\ss^{3+1}, g)$ such that there are two points in $\Sigma$  that belong to the same null geodesic that is $g$-orthogonal to $\Sigma.$

\begin{figure}[htbp]     
\begin{center}     
\epsfxsize 5cm     
\hepsffile{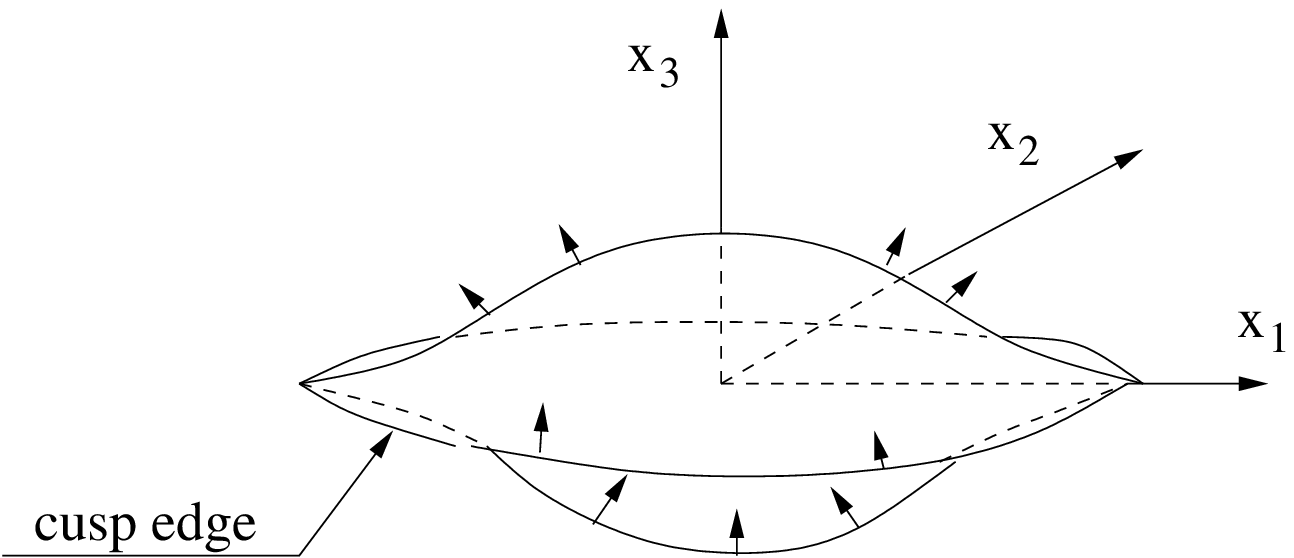}   
\end{center}     
\caption{}\label{nullexample.fig}     
\end{figure}

It also is possible to find Legendrian submanifolds of $\frak N$ that are not realizable as null congruences orthogonal to a spacelike $2$-surface. Consider 
a globally hyperbolic $\R^4$ with coordinates $(x_1, x_2, x_3, t)$ and the Lorentz metric  $g=dx_1^2+dx_2^2+dx_3^2-dt^2.$ For $\tau\in \R$ define the spacelike Cauchy surface $\R^3_{\tau}\subset \R^4$ to be the set of all the points whose $t$-coordinate equals $\tau.$ Take a Legendrian submanifold $L\subset ST^*\R^3_0=ST\R^3_0$ that is described by the projection of $L$ to $\R^3_0$ which is the rotationally symmetric ``flying saucer'' and the unit length vector field along the projection of $L$  orthogonal to the ``saucer'', see Figure~\ref{nullexample.fig}. We assume that the cusp edge of the ``saucer'' is the circle $\{(x_1, x_2, 0, 0)|x_1^2+x_2^2=1\}\subset \R^3_0\subset \R^4.$

Since $\R^3_0$ is a Cauchy surface, $ST^*\R^3_0$ is identified with $\frak N$ and $L$ gives a Legendrian submanifold $L'\subset \frak N.$ For every $\tau,$ the intersection of the corresponding null congruence with $\R^3_{\tau}$ will have a cusp edge along the circle 
$\{(x_1, x_2, {\tau}, {\tau})|x_1^2+x_2^2=1\}\subset \R^3_{\tau}\subset \R^4$ (and possibly other singularities). It is easy to see that the subsets of $L'$ that give rise to the cusp edges $\{(x_1, x_2, 0, 0)|x_1^2+x_2^2=1\}\subset \R^3_0$ and $\{(x_1, x_2, {\tau}, {\tau})|x_1^2+x_2^2=1\}\subset \R^3_{\tau}$ are equal. One verifies that no open in $L'$  neighborhood of a point in this subset can be realized as a null congruence orthogonal to some embedded (or immersed) $2$-dimensional spacelike surface in $(\R^4, g).$

Low also remarked that null congruences orthogonal to $2$-dimensional spacelike surfaces are related to Legendrian submanifolds of $\frak N$ for strongly causal $(\ss^{3+1}, g).$ In this case $\frak N$ is a smooth contact manifold that is possibly not Hausdorff. 

The contact structure on the space of null-geodesics and the symplectic structure on the spaces of timelike and spacelike geodesics in general pseudo-Riemannian manifolds was very recently studied by Khesin and Tabachnikov~\cite{KhesinTabachnikov}. In order for their results to apply the pseudo-Riemannian manifold should be such that these spaces of geodesics are manifolds. This imposes very strong restrictions on the pseudo-Riemannian manifolds under consideration.

Our work is motivated by Low's work and it establishes relations between Legendrian and null maps for an arbitrary space-time $(\ss^{m+1}, g),$ including those space-times for which the space of null-geodesics is not a manifold.
In particular, Theorem~\ref{mappednullsurface} shows that for an immersed spacelike hypersurface $\mu:\mathcal M^m\to \ss^{m+1},$ the null congruence associated to a Legendrian map $\wt \lambda:\mathcal L^{m-1}\to ST^*\mathcal M$ gives a mapped null hypersurface $\nu:\mathcal N^m\to \ss^{m+1},$ and that moreover the natural map $\wt \nu:\mathcal N^m\to ST^*\ss$ is a Legendrian immersion if $\wt \lambda$ is a Legendrian immersion.

Theorem~\ref{tripback} says that the intersection of a mapped null  hypersurface $\nu:\mathcal N^m\to \ss^{m+1}$ with any immersed spacelike hypersurface $\mu':\mathcal M'\to \ss$ gives a Legendrian map to $ST^*\mathcal M',$ and that moreover this map is a Legendrian immersion if $\wt \nu$ is an immersion. In this case the intersection $(\mu')^{-1}\bigl(\nu(\mathcal N)\cap \mu'(\mathcal M')\bigr)$ 
is naturally parameterized by a Legendrian immersion to $ST^*\mathcal M'.$

{\bf Acknowledgments:}  
I am very thankful to Robert Low, Yuli Rudyak, Miguel Sanchez, and David Webb for useful 
discussions.


\begin{thebibliography}{99999} 
 



\bibitem{BeemEhrlichEasley} 
J.~K.~Beem, P.~E.~Ehrlich, K.~L.~Easley: {\em Global Lorentzian geometry.\/} 
Second edition. Monographs and Textbooks in Pure and Applied Mathematics, {\bf 
202\/} Marcel Dekker, Inc., New York, 1996 

 
\bibitem{BernalSanchez}  
A.~Bernal, M.~Sanchez: {\it On smooth Cauchy hypersurfaces and Geroch's  
splitting  
theorem, \/} Comm. Math. Phys. {\bf 243\/} (2003), no. 3, 461--470 



%



\bibitem{Geroch}  
R.~P.~Geroch: {\it Domain of dependence,\/} J.~Math.~Phys., {\bf 11} (1970)  
pp.~437--449  
 

\bibitem{GromollKlingenbergMeyer} 
D.~Gromoll, W.~Klingenberg, W.~Meyer: {\em Riemannsche Geometrie im Grossen.\/} 
(German) Second edition, Lecture Notes in Mathematics, Vol. {\bf 55} 
Springer-Verlag, Berlin-New York (1975) 


\bibitem{HawkingEllis}
S.~W.~Hawking and G.~F.~R.~Ellis, {\em The large scale structure of space-time.\/} 
Cambridge Monographs on Mathematical Physics, No. 1. 
Cambridge University Press, London-New York, (1973)

\bibitem{KhesinTabachnikov}
B.~Khesin and S.~Tabachnikov: {\it Spaces of pseudo-Riemannian geodesics and pseudo-Euclidean billiards,\/} preprint math.DG/0608620 at http://www.arxiv.org, 46 pages (2006) 




\bibitem{Lerner} 
D.~Lerner: {\em Techniques of topology and differential geometry in general 
relativity.\/} 
Based on lectures by R. Penrose given at this conference. 
Methods of local and global differential geometry in general relativity (Proc. 
Regional Conf.~Relativity, Univ.~Pittsburgh, Pittsburgh, Pa., 1970), pp. 1--44. 
Lecture Notes in Physics, Vol. 14, 
Springer, Berlin (1972) 


\bibitem{Low0} 
R.~J.~Low: {\em Causal relations and spaces of null geodesics.\/} PhD Thesis, 
Oxford University, (1988) 


\bibitem{LowLegendrian} 
R.~J.~Low: {\em Stable singularities of wave-fronts in general relativity.\/} J. 
Math. Phys. 39 (1998), no. 6, 3332--3335 

\bibitem{NatarioTod} 
J.~Natario and P.~Tod: 
{\em Linking, Legendrian linking and causality.\/} 
Proc. London Math. Soc. (3) 88 (2004), no. 1, 251--272. 


\bibitem{ONeill} 
B.~O'Neill: {\it Semi-Riemannian geometry. With applications to relativity.\/} 
Pure and Applied Mathematics, 103. Academic Press, Inc. [Harcourt Brace 
Jovanovich, Publishers], New York, (1983) 

\bibitem{Whitehead}
J.H.C.~Whitehead: {\it Convex regions in the geometry of paths.\/}
Q. J. Math., Oxf. Ser. 3, 33-42 (1932)

\bibitem{Whiteheadaddendum}
J.H.C.~Whitehead: {\it Convex regions in the geometry of paths. Addendum.\/}
Q. J. Math., Oxf. Ser. 4, 226-227 (1933)
\end{thebibliography}
\end{document}